\documentclass{amsart}
\usepackage{amssymb}

\newtheorem{theorem}{Theorem}[section]
\newtheorem{proposition}[theorem]{Proposition}
\newtheorem{lemma}[theorem]{Lemma}

\theoremstyle{definition}
\newtheorem{definition}[theorem]{Definition}

\numberwithin{equation}{section}

\begin{document}
\title[$p$-Adic hypergeometric functions and twisted Kloosterman sheaf sum]
{$p$-Adic hypergeometric functions in the connections with certain twisted Kloosterman sheaf sum and modular forms}

{}
\author{Neelam Saikia}
\address{Department of Mathematics, Indian Institute of Technology Guwahati, North Guwahati, Guwahati-781039, Assam, INDIA}
\curraddr{}
\email{neelam16@iitg.ernet.in}
\thanks{}


\thanks{The author acknowledges
the financial support of Department of Science and Technology, Government of India for supporting a part of this work under INSPIRE Faculty Fellowship.}
\begin{abstract} 
In this paper we establish certain identities connecting $p$-adic hypergeometric functions with 4-th twisted Kloosterman sheaf sum. To prove these identities we express certain character sum over finite field in terms of special values of $p$-adic hypergeometric functions. One conjecture of Evans behaves as a bridge to connect $p$-adic hypergeometric functions with Kloosterman sheaf sum. We also connect  $p$-adic hypergeometric functions with Fourier coefficients of certain modular forms.
\end{abstract}
\maketitle
\section{Introduction and statement of results}
Let $p$ be an odd prime and $\mathbb{F}_p$ denote a finite field with $p$ elements. For any multiplicative character $\chi$ on $\mathbb{F}_p^{\times}$, we extend the domain of $\chi$ to a function on $\mathbb{F}_p$ by simply setting $\chi(0)=0$ including the trivial character $\varepsilon$. Let $\theta:\mathbb{F}_p\rightarrow\mathbb{C}$ be the additive character $x\mapsto e^{\frac{2\pi ix}{p}}$. 
\par Let $\overline{x}$ denote the multiplicative inverse of $x\in\mathbb{F}_p^{\times}$. 
For $a\in\mathbb{F}_p$, the classical Kloosterman sum is defined by
\begin{align}\label{eqn-1}
K(a):=\sum_{x\in\mathbb{F}_p^{\times}}\theta(x+a\overline{x}).
\end{align}
The $n$-th power Kloosterman moment is defined by $$S(n):=\displaystyle\sum_{a=0}^{p-1}K(a)^n.$$ For $a\neq0$, let $g(a)$ and $h(a)$ be the roots of the polynomial $X^2+K(a)X+p.$ Then one can write 
$$S(n)=(-1)^n+(-1)^n\displaystyle\sum_{a=1}^{p-1}(g(a)+h(a))^n.$$
The moments of the classical Kloosterman sums and their connection with Hecke eigenforms have been studied by many mathematicians, see for example \cite{ce, di, ronald}. For quadratic character $\varphi$ and for integer $n\geq1$, the $n$-th power twisted Kloosterman moment is defined as
\begin{align*}
S(n,\varphi):=\sum_{a\in\mathbb{F}_p^{\times}}\varphi(a)K(a)^n.
\end{align*}
To study the twisted Kloosterman moments, it is natural to consider the 
the $n$th twisted Kloosterman sheaf sum of $\varphi$ defined as follows:
\begin{align}\label{eqn-2}
T_{n,\varphi}:=\sum_{a\in\mathbb{F}_p^{\times}}\varphi(a)(g(a)^n+g(a)^{n-1}h(a)+g(a)^{n-2}h(a)^2+\cdots +g(a)h(a)^{n-1}+h(a)^n).
\end{align}
Let $\widehat{\mathbb{F}_p^\times}$ denote the group of all multiplicative characters on $\mathbb{F}_p^{\times}$. Let $\overline{\chi}$ denote the inverse of $\chi$. For multiplicative characters $A$ and $B$ on $\mathbb{F}_p^{\times}$ the Jacobi sum is defined by
\begin{align*}
J(A,B):=\sum_{x\in\mathbb{F}_p^{\times}}A(x)B(1-x)
\end{align*}
and the normalized Jacobi sum known as binomial is defined by
\begin{align}\label{binomial}
{A\choose B}:=\frac{B(-1)}{p}J(A,\overline{B}).
\end{align}
For multiplicative characters $A_0,A_1,\ldots, A_n$ and $B_1,B_2,\ldots,B_n$ on $\mathbb{F}_p^{\times}$ and 
$x\in\mathbb{F}_p$, Greene \cite{greene} first introduced the notion of hypergeometric function over finite fields by
\begin{align*}
{_{n+1}F_n}\left(\begin{array}{cccc}
A_0,& A_1,&\ldots,& A_n\\
~& B_1,&\ldots,& B_n
\end{array}|x
\right):=\frac{p}{p-1}\sum_{\chi\in\widehat{\mathbb{F}_p^{\times}}}{A_0\chi\choose \chi}
{A_1\chi\choose B_1\chi}\cdots{A_n\chi\choose B_n\chi}\chi(x).
\end{align*}
This function is also known as Gaussian hypergeometric function. Greene defined this function analogous to the classical hypergeometric series. We refer to \cite{greene} for a detailed study on Gaussian hypergeometric functions. There are other analogues of classical hypergeometric series over finite fields. For example, functions of this type were also defined by Katz \cite{katz} and McCarthy \cite{mccarthy3}. However, we can relate these finite field hypergeometric functions with each other, see for more details \cite{mccarthy3}. 
For brevity, we fix a notation for a special Gaussian hypergeometric function as follows:
\begin{align*}
{_{n+1}F_n}(x):={_{n+1}F_n}\left(\begin{array}{cccc}
\varphi,& \varphi,&\ldots,& \varphi\\
~& \varepsilon,&\ldots,& \varepsilon
\end{array}|x
\right).
\end{align*}
Gaussian hypergeometric functions play significant role in many branches of mathematics. For example, they appear in the theory of modular forms and elliptic curves \cite{ahl-ono-2, evans3, Fuselier, fuselier-thesis, fuselier2, fop, mccarthy4, mortenson}.
By definition, results involving hypergeometric functions over finite fields often be restricted to certain congruence classes of primes; see for example \cite{evans2, evans3, fuselier2, fuselier-thesis, mortenson}. To overcome these restrictions McCarthy \cite{mccarthy2, mccarthy1} defined a function $_nG_n[\cdots]$ in terms of Teichm\"{u}ller character over finite field and quotients of 
$p$-adic gamma functions which can be best described as an analogue of hypergeometric series in the $p$-adic setting. 
Let $\Gamma_p(\cdot)$ denote Morita's $p$-adic gamma function (defined in Section 2) and $\omega$ denote the Teichm\"{u}ller character of $\mathbb{F}_p$ with $\overline{\omega}$ denoting its character inverse. For $x\in\mathbb{Q}$ let $\lfloor x\rfloor$ denote the greatest integer less than or equal to $x$ and $\langle x\rangle$ denote the fractional part of $x$, satisfying $0\leq\langle x\rangle<1$.
We now recall the McCarthy's $p$-adic hypergeometric function $_{n}G_{n}[\cdots]$
as follows. 
\begin{definition}\cite[Definition 5.1]{mccarthy2} \label{defin1}
Let $p$ be an odd prime and $t \in \mathbb{F}_p$.
For positive integer $n$ and $1\leq k\leq n$, let $a_k$, $b_k$ $\in \mathbb{Q}\cap \mathbb{Z}_p$.
Then the function $_{n}G_{n}[\cdots]$ is defined as
\begin{align}
&_nG_n\left[\begin{array}{cccc}
             a_1, & a_2, & \ldots, & a_n \\
             b_1, & b_2, & \ldots, & b_n
           \end{array}|t
 \right]:=\frac{-1}{p-1}\sum_{a=0}^{p-2}(-1)^{an}~~\overline{\omega}^a(t)\notag\\
&\times \prod\limits_{k=1}^n(-p)^{-\lfloor \langle a_k \rangle-\frac{a}{p-1} \rfloor -\lfloor\langle -b_k \rangle +\frac{a}{p-1}\rfloor}
 \frac{\Gamma_p(\langle a_k-\frac{a}{p-1}\rangle)}{\Gamma_p(\langle a_k \rangle)}
 \frac{\Gamma_p(\langle -b_k+\frac{a}{p-1} \rangle)}{\Gamma_p(\langle -b_k \rangle)}.\notag
\end{align}
\end{definition}
This function is also known as $p$-adic hypergeometric function. 
One of the fundamental importance of $p$-adic hypergeometric functions is that using these functions often the results involving finite field hypergeometric functions can be extended to a wider class of primes, see for example \cite{BS1, mccarthy1, mccarthy2, mccarthy4}.
Let $q:=e^{2\pi iz}$, where $z\in\mathcal{H}$, the complex upper half plane.
Let $$f_1(z)=q\prod_{n=1}^{\infty}(1-q^{2n})^4(1-q^{4n})^4=\displaystyle\sum_{n=1}^{\infty}a(n)q^n$$ 
be the unique normalized newform in $S_4(\Gamma_0(8))$, the space of weight four cusp forms on the
congruence subgroup $\Gamma_0(8)$. Also, let $f_2\in S_4(\Gamma_0(16))$ be the unique weight 4 newform of level 16 with Fourier expansion $f_2(z)=\displaystyle\sum_{n=1}^{\infty}b(n)q^n$, which is the quadratic twist of $f_1$ with
$b(p)=\varphi(-1)a(p)$.\\\\
The $p$-adic hypergeometric functions is in the early state of development. 
These hypergeometric functions have appeared in different areas of mathematics. We observe that their connections with different objects are very significant. In particular, they have interesting applications in modular forms and elliptic curves. For example, Fuselier and McCarthy \cite{fm} expressed the Fourier coefficients of the modular form $f_2$ as special value of $p$-adic hypergeometric function as follows. 
\begin{align*}
{_4G_4}\left[\begin{array}{cccc}
             \frac{1}{2}, & \frac{1}{2}, & \frac{1}{4}, & \frac{3}{4} \\
             1, & 1, & 1, & 1
           \end{array}|1\right]-s(p)\cdot p=b(p),
\end{align*}
where $s(p)=\Gamma_p\left(\frac{1}{2}\right)^2\Gamma_p\left(\frac{1}{4}\right)\Gamma_p\left(\frac{3}{4}\right)$.
This identity was used to establish one of supercongruence conjecture of Rodriguez-Villegas between a classical truncated hypergeometric series and Fourier coefficients of the modular form $f_2$. Recently, Pujahari and the author \cite{SN} have  expressed the traces of $p$-th Hecke operators in terms of $p$-adic hypergeometric functions and applying this trace formula we connected Ramanujan's $\tau$ function with $p$-adic hypergeometric functions.
These consequences motivates the author to explore more connections between $p$-adic hypergeometric functions and Fourier coefficients of modular forms.\\\\
Many authors studied the connections between Kloosterman sums and hypergeometric functions and established many motivating results. In \cite{evans3} Evans gave relations connecting certain moments of the twisted Kloosterman sums with Gaussian hypergeometric functions and gave some conjectures between the twisted Kloosterman sums and the Fourier coefficients of modular forms. Dummit, Goldberg, and Perry \cite{DGP} established the 4-th quadratic twisted Kloosterman sheaf sum as a special value of the Gaussian hypergeometric function ${_4F_3}(1)$. Then combining their results with results of Ahlgren and Ono \cite{ahl-ono-1, ahl-ono-2}, Dummit, Goldberg, and Perry \cite{DGP} proved one conjecture of of Evans.
Applying \cite[Theorem 1.1]{DGP}, McCarthy's transformations \cite[Lemma 3.3]{mccarthy2} and \cite[Proposition 2.5]{mccarthy3} we can write $p$-adic hypergeometric function in terms of 4-th quadratic twisted Kloosterman sheaf sum as follows:
\begin{align}
p^6~{_4G_4}\left[\begin{array}{cccc}
\frac{1}{2}, & \frac{1}{2}, & \frac{1}{2}, & \frac{1}{2}\\
0, & 0, & 0, & 0
\end{array}\mid 1\right]=p-\frac{T_{4,\varphi}}{p}.\notag
\end{align}
Also, Lin and Tu \cite{lin-tu} established certain identities among twisted moments of twisted Kloosterman sums, and hypergeometric functions over finite fields. It seems that connections of hypergeometric functions with different objects have been interesting and beneficial in the literature. All these connections motives us to explore more relations between $p$-adic hypergeometric functions and Kloosterman sheaf sums.\\\\
The main goal of the paper is to establish certain identities expressing special values of $p$-adic hypergeometric functions in terms of 4-th quadratic twisted Kloosterman sheaf sum $T_{4,\varphi}$.  These identities unable us to connect $p$-adic hypergeometric functions with the $p$-th Fourier coefficients of modular forms $f_1$ and $f_2$.\\\\
We now state our first result which provides two relations connecting a special value of ${_4G_4}[\cdots]$ function with $T_{4,\varphi}$ under the condition that $p\equiv1$ or $7\pmod{12}$. Let $\psi_6=\omega^{\frac{p-1}{6}}$ and 
$\psi_3=\omega^{\frac{p-1}{3}}$ be characters of order $6$ and $3$, respectively.
\begin{theorem}\label{MT-1}
If $p\equiv1\pmod{12}$ such that $p=x^2+y^2$ and $x$ is odd, then we have
\begin{align}\label{case-1}
&Cp^3(p-1)\overline{\psi}_6(2)\psi_3(4){_4G_4}\left[\begin{array}{cccc}
             \frac{5}{6}, & \frac{5}{6}, & \frac{1}{12}, & \frac{7}{12}\vspace{.1 cm } \\
             \frac{1}{3}, & \frac{1}{3}, & \frac{1}{3}, & \frac{1}{3}
           \end{array}|1\right]\notag\\
           &=1-T_{4,\varphi}-p^2\varphi(2)+4px^2\varphi(2)-p\varphi(-2)(4x^2-2p),
\end{align}
and if $p\equiv7\pmod{12}$, then we have 
\begin{align}\label{case-2}
Cp^3(p-1)\overline{\psi}_6(2)\psi_3(4)
{_4G_4}\left[\begin{array}{cccc}
             \frac{5}{6}, & \frac{5}{6}, & \frac{1}{12}, & \frac{7}{12}\vspace{.1 cm } \\
             \frac{1}{3}, & \frac{1}{3}, & \frac{1}{3}, & \frac{1}{3}
           \end{array}|1\right]&=1-T_{4,\varphi}-p^2\varphi(2),
\end{align}
where $C=\frac{\Gamma_p(\frac{2}{3})^4\Gamma_p(\frac{5}{6})^2\Gamma_p(\frac{1}{12})
\Gamma_p(\frac{7}{12})}{\Gamma_p(\frac{1}{2})}$.
\end{theorem}
The following theorem we link a special value of ${_{12}G_{12}}[\cdots]$-function with $T_{4,\varphi}$ for primes $p$ which are missed in Theorem \ref{MT-1}.
\begin{theorem}\label{MT-2}
If $p\equiv5\pmod{12}$ such that $p=x^2+y^2$ and $x$ is odd, then we have 
\begin{align}\label{case-3}
&{_{12}G_{12}}\left[\begin{array}{cccccccccccc}
             0, & \frac{1}{3}, & \frac{2}{3}, & 0, & \frac{1}{3}, & \frac{2}{3}, & 0, & \frac{1}{3}, & \frac{2}{3}, &
             0, & \frac{1}{3}, & \frac{2}{3} \vspace{.1 cm } \\
             \frac{1}{6}, & \frac{1}{2}, & \frac{5}{6}, & \frac{1}{12}, & \frac{1}{6}, & \frac{1}{4},
             & \frac{5}{12}, & \frac{1}{2}, & \frac{7}{12}, & \frac{3}{4}, & \frac{5}{6}, & \frac{11}{12}
           \end{array}\mid1\right]\notag\\
           &=\frac{1}{p-1}
\left(\varphi(-1)-\varphi(-1)T_{4,\varphi}-
p^2\varphi(-2)+4px^2\varphi(-2)
-p\varphi(2)(4x^2-2p)\right),
\end{align}
and if $p\equiv11\pmod{12}$, then we have
\begin{align}\label{case-4}
&{_{12}G_{12}}\left[\begin{array}{cccccccccccc}
             0, & \frac{1}{3}, & \frac{2}{3}, & 0, & \frac{1}{3}, & \frac{2}{3}, & 0, & \frac{1}{3}, & \frac{2}{3}, &
             0, & \frac{1}{3}, & \frac{2}{3} \vspace{.1 cm } \\
             \frac{1}{6}, & \frac{1}{2}, & \frac{5}{6}, & \frac{1}{12}, & \frac{1}{6}, & \frac{1}{4},
             & \frac{5}{12}, & \frac{1}{2}, & \frac{7}{12}, & \frac{3}{4}, & \frac{5}{6}, & \frac{11}{12}
           \end{array}\mid1\right]\notag\\
&=\frac{1}{p-1}\left(\varphi(-1)-\varphi(-1)T_{4,\varphi}-p^2\varphi(-2)\right).
\end{align}
\end{theorem}
As a consequence of Theorem \ref{MT-1} we obtain the following theorem, which expresses the  $p$-th Fourier coefficients of $f_1$ and $f_2$ in terms of $p$-adic hypergeometric functions.
\begin{theorem}\label{MT-3}
If $p\equiv1\pmod{12}$ such that $p=x^2+y^2$ and $x$ is odd, then we have
\begin{align}\label{a-1}
a(p)&=-\frac{1}{p}+p\varphi(2)-4x^2\varphi(2)+\varphi(-2)(4x^2-2p)\notag\\
&~~+Cp^2(p-1)\overline{\psi}_6(2)\psi_3(4){_4G_4}\left[\begin{array}{cccc}
             \frac{5}{6}, & \frac{5}{6}, & \frac{1}{12}, & \frac{7}{12}\vspace{.1 cm } \\
             \frac{1}{3}, & \frac{1}{3}, & \frac{1}{3}, & \frac{1}{3}
           \end{array}|1\right],
\end{align}
and 
\begin{align}\label{b-1}
b(p)&=-\frac{\varphi(-1)}{p}+p\varphi(-2)-4x^2\varphi(-2)+\varphi(2)(4x^2-2p)\notag\\
&~~+Cp^2(p-1)\varphi(-1)\overline{\psi}_6(2)\psi_3(4){_4G_4}\left[\begin{array}{cccc}
             \frac{5}{6}, & \frac{5}{6}, & \frac{1}{12}, & \frac{7}{12}\vspace{.1 cm } \\
             \frac{1}{3}, & \frac{1}{3}, & \frac{1}{3}, & \frac{1}{3}
           \end{array}|1\right].
\end{align}
If $p\equiv7\pmod{12}$, then we have 
\begin{align}\label{a-2}
a(p)&=-\frac{1}{p}+p\varphi(2)+Cp^2(p-1)\overline{\psi}_6(2)\psi_3(4){_4G_4}\left[\begin{array}{cccc}
             \frac{5}{6}, & \frac{5}{6}, & \frac{1}{12}, & \frac{7}{12}\vspace{.1 cm } \\
             \frac{1}{3}, & \frac{1}{3}, & \frac{1}{3}, & \frac{1}{3}
           \end{array}|1\right],
\end{align}
and 
\begin{align}\label{b-2}
b(p)&=-\frac{\varphi(-1)}{p}+p\varphi(-2)+Cp^2(p-1)\varphi(-1)\overline{\psi}_6(2)\psi_3(4){_4G_4}\left[\begin{array}{cccc}
             \frac{5}{6}, & \frac{5}{6}, & \frac{1}{12}, & \frac{7}{12}\vspace{.1 cm } \\
             \frac{1}{3}, & \frac{1}{3}, & \frac{1}{3}, & \frac{1}{3}
           \end{array}|1\right].
\end{align}
\end{theorem}
We also prove the following theorem as a consequence of Theorem \ref{MT-2}. The theorem connects the $p$-th Fourier coefficients of $f_1$ and $f_2$ in terms of special values of $p$-adic hypergeometric function for primes which are not included in Theorem \ref{MT-3}.
\begin{theorem}\label{MT-4}
If $p\equiv5\pmod{12}$ such that $p=x^2+y^2$ and $x$ is odd, then we have 
\begin{align}\label{a-3}
a(p)&=-\frac{1}{p}+p\varphi(2)-4x^2\varphi(2)+\varphi(-2)(4x^2-2p)+\frac{\varphi(-1)(p-1)}{p}\notag\\
&~~\times{_{12}G_{12}}\left[\begin{array}{cccccccccccc}
             0, & \frac{1}{3}, & \frac{2}{3}, & 0, & \frac{1}{3}, & \frac{2}{3}, & 0, & \frac{1}{3}, & \frac{2}{3}, &
             0, & \frac{1}{3}, & \frac{2}{3} \vspace{.1 cm } \\
             \frac{1}{6}, & \frac{1}{2}, & \frac{5}{6}, & \frac{1}{12}, & \frac{1}{6}, & \frac{1}{4},
             & \frac{5}{12}, & \frac{1}{2}, & \frac{7}{12}, & \frac{3}{4}, & \frac{5}{6}, & \frac{11}{12}
           \end{array}\mid1\right],
\end{align}
and 
\begin{align}\label{b-3}
b(p)&=-\frac{\varphi(-1)}{p}+p\varphi(-2)-4x^2\varphi(-2)+\varphi(2)(4x^2-2p)
\notag\\
&+\frac{p-1}{p}~{_{12}G_{12}}\left[\begin{array}{cccccccccccc}
             0, & \frac{1}{3}, & \frac{2}{3}, & 0, & \frac{1}{3}, & \frac{2}{3}, & 0, & \frac{1}{3}, & \frac{2}{3}, &
             0, & \frac{1}{3}, & \frac{2}{3} \vspace{.1 cm } \\
             \frac{1}{6}, & \frac{1}{2}, & \frac{5}{6}, & \frac{1}{12}, & \frac{1}{6}, & \frac{1}{4},
             & \frac{5}{12}, & \frac{1}{2}, & \frac{7}{12}, & \frac{3}{4}, & \frac{5}{6}, & \frac{11}{12}
           \end{array}\mid1\right].
\end{align}
If $p\equiv11\pmod{12}$, then we have
\begin{align}\label{a-4}
a(p)&=-\frac{1}{p}+p\varphi(2)\notag\\
&~~\times{_{12}G_{12}}\left[\begin{array}{cccccccccccc}
             0, & \frac{1}{3}, & \frac{2}{3}, & 0, & \frac{1}{3}, & \frac{2}{3}, & 0, & \frac{1}{3}, & \frac{2}{3}, &
             0, & \frac{1}{3}, & \frac{2}{3} \vspace{.1 cm } \\
             \frac{1}{6}, & \frac{1}{2}, & \frac{5}{6}, & \frac{1}{12}, & \frac{1}{6}, & \frac{1}{4},
             & \frac{5}{12}, & \frac{1}{2}, & \frac{7}{12}, & \frac{3}{4}, & \frac{5}{6}, & \frac{11}{12}
           \end{array}\mid1\right],
\end{align}
and
\begin{align}\label{b-4}
b(p)&=-\frac{\varphi(-1)}{p}+p\varphi(-2)\notag\\
&+\frac{p-1}{p}~{_{12}G_{12}}\left[\begin{array}{cccccccccccc}
             0, & \frac{1}{3}, & \frac{2}{3}, & 0, & \frac{1}{3}, & \frac{2}{3}, & 0, & \frac{1}{3}, & \frac{2}{3}, &
             0, & \frac{1}{3}, & \frac{2}{3} \vspace{.1 cm } \\
             \frac{1}{6}, & \frac{1}{2}, & \frac{5}{6}, & \frac{1}{12}, & \frac{1}{6}, & \frac{1}{4},
             & \frac{5}{12}, & \frac{1}{2}, & \frac{7}{12}, & \frac{3}{4}, & \frac{5}{6}, & \frac{11}{12}
           \end{array}\mid1\right].
\end{align}
\end{theorem}
We now briefly outline the structure of the paper. In Section 2 we give necessary definitions and preliminary results. We state certain results including Davenport-Hasse relation and Gross-Koblitz formula in Section 2. We also prove certain lemmas and  propositions in Section 2. Finally in Section 3 we prove Theorem \ref{MT-1}--\ref{MT-4}. 
\section{Notation and Preliminaries}
\subsection{Multiplicative character, Gauss sum and Jacobi sum}
In this section we start with a lemma that provides the orthogonality relations satisfied by multiplicative characters.
\begin{lemma}\emph{(\cite[Chapter 8]{ireland}).}\label{lemma-3} We have
\begin{enumerate}
\item $\displaystyle\sum\limits_{x\in\mathbb{F}_p}\chi(x)=\left\{
                                  \begin{array}{ll}
                                    p-1 & \hbox{if~ $\chi=\varepsilon$;} \\
                                    0 & \hbox{if ~~$\chi\neq\varepsilon$.}
                                  \end{array}
                                \right.$
\item $\displaystyle\sum\limits_{\chi\in \widehat{\mathbb{F}_p^\times}}\chi(x)~~=\left\{
                            \begin{array}{ll}
                              p-1 & \hbox{if~~ $x=1$;} \\
                              0 & \hbox{if ~~$x\neq1$.}
                            \end{array}
                          \right.$
\end{enumerate}
\end{lemma}
\par
Let $\mathbb{Z}_p$ and $\mathbb{Q}_p$ denote the ring of $p$-adic integers and the field of $p$-adic numbers, respectively.
Let $\overline{\mathbb{Q}_p}$ be the algebraic closure of $\mathbb{Q}_p$ and $\mathbb{C}_p$ be the completion of $\overline{\mathbb{Q}_p}$. We now introduce Gauss sum and some interesting properties of Gauss sums. For further details, see \cite{evans}.
Let $\zeta_p$ be a fixed primitive $p$-th root of unity in $\overline{\mathbb{Q}_p}$. Then the additive character
$\theta: \mathbb{F}_p \rightarrow \mathbb{Q}_p(\zeta_p)$ is defined by
\begin{align}
\theta(\alpha)=\zeta_p^{\alpha}.\notag
\end{align}
It is easy to check that the additive character $\theta$ satisfies the following relations:
$$\theta(a+b)=\theta(a)\theta(b)$$ and
\begin{align}\label{lemma-4}
\sum_{x\in\mathbb{F}_p}\theta(x)=0.
\end{align}
For $\chi \in \widehat{\mathbb{F}_p^\times}$, the Gauss sum is defined by
\begin{align}
g(\chi):=\sum\limits_{x\in \mathbb{F}_p}\chi(x)\theta(x).\notag
\end{align}
From \eqref{lemma-4} we can say that $g(\varepsilon)=-1$.
The following elementary property of Gauss sum is very useful, which gives a nice expression for the product of two Gauss sums.
\begin{lemma}\emph{(\cite[Eq. 1.12]{greene}).}\label{lemma-1}
Let $\chi\in\widehat{\mathbb{F}_p^{\times}}$ be such that $\chi\neq\varepsilon$. Then
$$g(\chi)g(\overline{\chi})=p\cdot\chi(-1).$$
\end{lemma}
Also, another interesting product of Gauss sums is due to Davenport and Hasse. This relation is very well known as Davenport-Hasse formula.
\begin{theorem}\emph{(\cite[Davenport-Hasse Relation]{evans}).}\label{DH}
Let $m$ be a positive integer and let $p$ be a prime such that $p\equiv 1 \pmod{m}$. For multiplicative characters
$\chi, \psi \in \widehat{\mathbb{F}_p^\times}$, we have
\begin{align}
\prod\limits_{\chi^m=\varepsilon}g(\chi \psi)=-g(\psi^m)\psi(m^{-m})\prod\limits_{\chi^m=\varepsilon}g(\chi).\notag
\end{align}
\end{theorem}
Let $A,B$ be two multiplicative characters on $\mathbb{F}_p^{\times}$ and let $a\in\mathbb{F}_p^{\times}$. Then the  generalized Jacobi sum $J_a(A,B)$ is defined as
$$J_a(A,B):=\displaystyle\sum_{t_1+t_2=a}A(t_1)B(t_2).$$
If we put $t_1=ab_1$ and $t_2=ab_2$. Then $b_1+b_2=1$ and 
\begin{align}\label{eqn-J}
J_a(A,B)&=\sum_{b_1+b_2=1}A(ab_1)B(ab_2)\notag\\
&=A(a)B(a)\sum_{b_1+b_2=1}A(b_1)B(b_2)\notag\\
&=AB(a)J(A,B).
\end{align}
We need the following elementary properties from \cite{greene} to prove our main results: Let $\delta$ be defined by
\begin{align*}
\delta(A)=\left\{
   \begin{array}{ll}
     1, & \hbox{if $A=\varepsilon$;} \\
    0, & \hbox{otherwise.}
   \end{array}
 \right.
\end{align*}
\begin{align}\label{property-1}
{A\choose B}&=B(-1){B\overline{A}\choose B},\\
\label{property-3}
{\varphi\choose \varphi}&=\frac{-1}{p}.
\end{align}
The following relation relates binomial with Gauss sums.
\begin{align}\label{property-2}
{A\choose B}&=\frac{g(A)g(\overline{B})B(-1)}{g(A\overline{B})p}+\frac{p-1}{p}\delta(A\overline{B}).
\end{align}
For more details on the properties of binomilas we refer \cite{greene}.

\subsection{Gaussian hypergeometric functions}
Here we state some basic results on Gaussian hypergeometric functions, for more details see \cite{greene}. We recall that for $x\neq0$,
\begin{align}
{_2F_1}(x)&=\frac{\varphi(-1)}{p}\sum_{y\in\mathbb{F}_p^{\times}}\varphi(y)\varphi(1-y)\varphi(1-xy),\\
{_3F_2}(x)&=\frac{1}{p^2}\sum_{y,z\in\mathbb{F}_p^{\times}}\varphi(y)\varphi(1-y)\varphi(z)\varphi(1-z)
\varphi(1-xyz).
\end{align}
Gaussian hypergeometric functions satisfy many transformation identities analogous to classical hypergeometric series. For instance, we consider the following transformation due to Evans and Greene, which is a special case of \cite[Theorem 1.7]{EG}. 
\begin{lemma}\label{lemma 10}
If $t\neq0,\pm1$, then we have
\begin{align}
{_3F_2}\left(\frac{1}{1-t^2}\right)=\varphi(t^2-1)\left(-\frac{1}{p}+{_2F_1}\left(\frac{1-t}{2}\right)^2\right).\notag
\end{align}
\end{lemma}
Recall two special values of the functions ${_2F_1}(x)$ and ${_3F_2}(x)$ from \cite{ono}. These values are very useful to prove our main results.
\begin{proposition}\label{value-1}\cite[Theorem 2]{ono}
Let $\lambda\in\{-1,\frac{1}{2},2\}$. If $p$ is an odd prime, then
\begin{align}
{_2F_1}\left(\lambda\right)=\left\{
   \begin{array}{ll}
     0, & \hbox{if $p\equiv3\pmod{4}$;} \\
    \frac{2x(-1)^{\frac{x+y+1}{2}}}{p}, & \hbox{if $p\equiv1\pmod{4}, x^2+y^2=p$ and $x$ odd.}
   \end{array}
 \right.\notag
\end{align}
\end{proposition}
\begin{proposition}\label{value-2}\cite[Theorem 4]{ono}
Let $p$ be an odd prime. Then 
\begin{align}
{_3F_2}(1)=\left\{
   \begin{array}{ll}
     0, & \hbox{if $p\equiv3\pmod{4}$;} \\
    \frac{4x^2-2p}{p^2}, & \hbox{if $p\equiv1\pmod{4}, x^2+y^2=p$ and $x$ odd.}
   \end{array}
 \right.\notag
\end{align}
\end{proposition}
We also use the following special value of Gaussian hypergeometric function from \cite{greene}.
\begin{proposition}\label{value-3}
Let $p$ be an odd prime. Then we have
\begin{align}
{_2F_1}(1)=-\frac{\varphi(-1)}{p}.\notag
\end{align}
\end{proposition}
\begin{proof}
From \cite[Theorem 4.9]{greene} it is well known that
${_2F_1(1)}=\varphi(-1){\varphi\choose\varphi}$. Then using \eqref{property-3} we obtain the result.
\end{proof}
\subsection{$p$-adic hypergeometric function}
We recall the $p$-adic gamma function, for further details see \cite{kob}.
For a positive integer $n$,
the $p$-adic gamma function $\Gamma_p(n)$ is defined as
\begin{align}
\Gamma_p(n):=(-1)^n\prod\limits_{0<j<n,p\nmid j}j\notag
\end{align}
and one extends it to all $x\in\mathbb{Z}_p$ by setting $\Gamma_p(0):=1$ and
\begin{align}
\Gamma_p(x):=\lim_{x_n\rightarrow x}\Gamma_p(x_n)\notag
\end{align}
for $x\neq0$, where $x_n$ runs through any sequence of positive integers $p$-adically approaching $x$. We now recall certain basic properties of $p$-adic gamma function. For example,
if $x\in\mathbb{Z}_p$, then we have
\begin{align}\label{prop-1}
\Gamma_p(1-x)\Gamma_p(x)=(-1)^{a_0(x)},
\end{align}
where $a_0(x)\in\{1,2,\ldots,p\}$ such that $x\equiv a_0(x)\pmod p$. The following product formula is used in the proof of Proposition \ref{prop-mt-1}. Let $m$ be a positive integer such that $p\nmid m$ and $x=\frac{r}{p-1}$ with $0\leq r\leq p-1$. Then
\begin{align}\label{product-1}
\prod_{h=0}^{m-1}\Gamma_p\left(\frac{x+h}{m}\right)=\omega^{(1-p)(1-x)}(m)\Gamma_p(x)\prod_{h=1}^{m-1}\Gamma_p\left(\frac{h}{m}\right).
\end{align}
For $x \in \mathbb{Q}$ we let $\lfloor x\rfloor$ denote the greatest integer less than
or equal to $x$ and $\langle x\rangle$ denote the fractional part of $x$, i.e., $x-\lfloor x\rfloor$, satisfying $0\leq\langle x\rangle<1$.
We now recall the McCarthy's $p$-adic hypergeometric function $_{n}G_{n}[\cdots]$
as follows.
\begin{definition}\cite[Definition 5.1]{mccarthy2} \label{defin1}
Let $p$ be an odd prime and let $t \in \mathbb{F}_p$.
For positive integer $n$ and $1\leq k\leq n$, let $a_k$, $b_k$ $\in \mathbb{Q}\cap \mathbb{Z}_p$.
Then the function $_{n}G_{n}[\cdots]$ is defined by
\begin{align}
&_nG_n\left[\begin{array}{cccc}
             a_1, & a_2, & \ldots, & a_n \\
             b_1, & b_2, & \ldots, & b_n
           \end{array}|t
 \right]:=\frac{-1}{p-1}\sum_{a=0}^{p-2}(-1)^{an}~~\overline{\omega}^a(t)\notag\\
&\times \prod\limits_{k=1}^n(-p)^{-\lfloor \langle a_k \rangle-\frac{a}{p-1} \rfloor -\lfloor\langle -b_k \rangle +\frac{a}{p-1}\rfloor}
 \frac{\Gamma_p(\langle a_k-\frac{a}{p-1}\rangle)}{\Gamma_p(\langle a_k \rangle)}
 \frac{\Gamma_p(\langle -b_k+\frac{a}{p-1} \rangle)}{\Gamma_p(\langle -b_k \rangle)}.\notag
\end{align}
\end{definition}
Let $\pi \in \mathbb{C}_p$ be the fixed root of the polynomial $x^{p-1} + p$, which satisfies the congruence condition
$\pi \equiv \zeta_p-1 \pmod{(\zeta_p-1)^2}$. Then the Gross-Koblitz formula allows us to relate Gauss sum with $p$-adic gamma function as follows.
\begin{theorem}\emph{(\cite[Gross-Koblitz]{gross}).}\label{GK} For $a\in \mathbb{Z}$,
\begin{align}
g(\overline{\omega}^a)=-\pi^{(p-1)\langle\frac{a}{p-1} \rangle}\Gamma_p\left(\left\langle \frac{a}{p-1} \right\rangle\right).\notag
\end{align}
\end{theorem}
The following lemma relates products of values of $p$-adic gamma function. This lemma is used in the proof of Proposition \ref{prop-mt-2}.
\begin{lemma}\emph{(\cite[Lemma 4.1]{mccarthy2}).}\label{lemma-5}
Let $p$ be a prime and $0\leq a\leq p-2$. For $t\geq 1$ with $p\nmid t$, we have
\begin{align}
\omega(t^{ta})\Gamma_p\left(\left\langle \frac{ta}{p-1}\right\rangle\right)
\prod\limits_{h=1}^{t-1}\Gamma_p\left(\left\langle\frac{h}{t}\right\rangle\right)
=\prod\limits_{h=0}^{t-1}\Gamma_p\left(\left\langle\frac{h}{t}+\frac{a}{p-1}\right\rangle\right),\notag
\end{align}
and
\begin{align}
\omega(t^{-ta})\Gamma_p\left(\left\langle\frac{-ta}{p-1}\right\rangle\right)
\prod\limits_{h=1}^{t-1}\Gamma_p\left(\left\langle \frac{h}{t}\right\rangle\right)
=\prod\limits_{h=1}^{t}\Gamma_p\left(\left\langle\frac{h}{t}-\frac{a}{p-1}\right\rangle \right).\notag
\end{align}
\end{lemma}
The next two elementary lemmas are very useful in the proof of Proposition \ref{prop-mt-2}.
\begin{lemma}\label{D1}
Let $p$ be an odd prime. Then for $d>0$ and $0\leq a\leq p-2$ we have 
\begin{align}
\left\lfloor\frac{-da}{p-1}\right\rfloor=-1+\sum_{h=1}^{d}\left\lfloor\left\langle\frac{h}{d}\right\rangle
-\frac{a}{p-1}\right\rfloor.
\end{align}
\end{lemma}
\begin{proof}
Let $\left\lfloor\frac{-da}{p-1}\right\rfloor=dl+k$ for some $l,k\in\mathbb{Z}$ such that $0\leq k\leq d-1$. Then we have the following inequality
\begin{align}\label{req-1}
l+\frac{k}{d}\leq-\frac{a}{p-1}<l+\frac{k+1}{d}.
\end{align}
The above inequality gives the following:
\begin{align}\label{req-2}
l+\frac{k+h}{d}\leq\frac{h}{d}-\frac{a}{p-1}<l+\frac{k+h+1}{d}
\end{align}
for $1\leq h\leq d$. Then we have
\begin{align}
\left\lfloor\left\langle\frac{h}{d}\right\rangle
-\frac{a}{p-1}\right\rfloor=\left\{
   \begin{array}{ll}
     l, & \hbox{if $h=1,2,\ldots,d-k-1$;} \\
    l+1, & \hbox{if $h=d-k,\ldots,d$,}
   \end{array}
 \right.
\end{align}
and this implies that 
$\displaystyle\sum_{h=1}^{d}\left\lfloor\left\langle\frac{h}{d}\right\rangle
-\frac{a}{p-1}\right\rfloor=dl+k+1=1+\left\lfloor\frac{-da}{p-1}\right\rfloor.$ This completes the proof of the lemma.
\end{proof}
\begin{lemma}\label{D2}
Let $p$ be an odd prime. Then for $d>0$ and $0\leq a\leq p-2$ we have 
\begin{align}
\left\lfloor\frac{da}{p-1}\right\rfloor=\sum_{h=0}^{d-1}\left\lfloor\left\langle\frac{h}{d}\right\rangle
+\frac{a}{p-1}\right\rfloor.
\end{align}
\end{lemma}
\begin{proof}
Let $\left\lfloor\frac{da}{p-1}\right\rfloor=dl+k$ for some $l,k\in\mathbb{Z}$ such that $0\leq k\leq d-1$. Then
following similar steps as in the proof of Lemma \ref{D1} we easily prove the lemma.
\end{proof}
The next two propositions play an important role in the proof of Theorem \ref{MT-1} and Theorem \ref{MT-2}. These two propositions express certain  sum of products of Gauss sums as special values of $p$-adic hypergeometric functions.
\begin{proposition}\label{prop-mt-1}
Let $p$ be an odd prime such that $p\equiv1\pmod{3}$. Then we have
\begin{align*}
&\sum_{a=0}^{p-2}\frac{g(\varphi\omega^{a})^2g(\overline{\omega}^a)^4g(\varphi\omega^{2a})}{g(\varphi)}~\overline{\omega}^a(4)\\&=-Cp^4(p-1)\overline{\psi}_6(2)\psi_3(4)~
{_4G_4}\left[\begin{array}{cccc}
             \frac{5}{6}, & \frac{5}{6}, & \frac{1}{12}, & \frac{7}{12}\vspace{.1 cm } \\
             \frac{1}{3}, & \frac{1}{3}, & \frac{1}{3}, & \frac{1}{3}
           \end{array}|1\right],
\end{align*}
where $C=\frac{\Gamma_p(\frac{2}{3})^4\Gamma_p(\frac{5}{6})^2\Gamma_p(\frac{1}{12})
\Gamma_p(\frac{7}{12})}{\Gamma_p(\frac{1}{2})}$ and where $\psi_6=\omega^{\frac{p-1}{6}}$ and 
$\psi_3=\omega^{\frac{p-1}{3}}$ are characters of order $6$ and $3$, respectively.
\end{proposition}
\begin{proof}
Let $B'=\displaystyle\sum_{a=0}^{p-2}\frac{g(\varphi\omega^{a})^2g(\overline{\omega}^a)^4g(\varphi\omega^{2a})}{g(\varphi)}~\overline{\omega}^a(4)$.
Then replacing $a$ by $a-\frac{p-1}{3}$ we obtain 
\begin{align*}
B'=\sum_{a=0}^{p-2}\frac{g(\varphi\omega^{a-\frac{p-1}{3}})^2g(\overline{\omega}^{a-\frac{p-1}{3}})^4g(\varphi\omega^{2a-\frac{2(p-1)}{3}})}{g(\varphi)}~\overline{\omega}^{a-\frac{p-1}{3}}(4).
\end{align*}
Applying Gross-Koblitz formula on each of the Gauss sums present in the above sum we have
\begin{align}\label{eqn-prop-1}
B&'=\frac{\psi_3(4)}{\Gamma_p(\frac{1}{2})}\sum_{a=0}^{p-2}\pi^{(p-1)S}~\overline{\omega}^{a}(4)
~\Gamma_p\left(\left\langle\frac{1}{2}+\frac{1}{3}-\frac{a}{p-1}\right\rangle\right)^2\Gamma_p\left(\left\langle\frac{2}{3}+\frac{a}{p-1}\right\rangle\right)^4\notag\\
&\times\Gamma_p\left(\left\langle\frac{1}{2}+\frac{2}{3}-\frac{2a}{p-1}\right\rangle\right)\notag\\
&=\frac{\psi_3(4)}{\Gamma_p(\frac{1}{2})}\sum_{a=0}^{p-2}\pi^{(p-1)S}~\overline{\omega}^{a}(4)
~\Gamma_p\left(\left\langle\frac{5}{6}-\frac{a}{p-1}\right\rangle\right)^2\Gamma_p\left(\left\langle\frac{2}{3}+\frac{a}{p-1}\right\rangle\right)^4\notag\\
&\times\Gamma_p\left(\left\langle\frac{1}{6}-\frac{2a}{p-1}\right\rangle\right),
\end{align}
where $S=2\left\langle\frac{5}{6}-\frac{a}{p-1}\right\rangle+4\left\langle\frac{2}{3}+\frac{a}{p-1}\right\rangle
+\left\langle\frac{1}{6}-\frac{2a}{p-1}\right\rangle-\frac{1}{2}$. If we apply \eqref{product-1} for $m=2$ and 
$x=\left\langle\frac{1}{6}-\frac{2a}{p-1}\right\rangle$, then we obtain
$$\Gamma_p\left(\left\langle\frac{1}{6}-\frac{2a}{p-1}\right\rangle\right)
=\frac{\omega^a(4)\overline{\psi}_6(2)}{\Gamma_p(\frac{1}{2})}
\Gamma_p\left(\left\langle\frac{1}{12}-\frac{a}{p-1}\right\rangle\right)
\Gamma_p\left(\left\langle\frac{7}{12}-\frac{a}{p-1}\right\rangle\right).$$ Substituting the above expression into 
\eqref{eqn-prop-1} we obtain
\begin{align}\label{eqn-prop-2}
B'&=\frac{\psi_3(4)\overline{\psi}_6(2)}{\Gamma_p(\frac{1}{2})^2}\sum_{a=0}^{p-2}\pi^{(p-1)S}
~\Gamma_p\left(\left\langle\frac{5}{6}-\frac{a}{p-1}\right\rangle\right)^2\Gamma_p\left(\left\langle\frac{2}{3}+\frac{a}{p-1}\right\rangle\right)^4\notag\\
&\times\Gamma_p\left(\left\langle\frac{1}{12}-\frac{a}{p-1}\right\rangle\right)
\Gamma_p\left(\left\langle\frac{7}{12}-\frac{a}{p-1}\right\rangle\right).
\end{align}
If we put $x=\frac{1}{2}$ into \eqref{prop-1} we obtain $\Gamma_p\left(\frac{1}{2}\right)^2=(-1)^{\frac{(p+1)}{2}}=-\varphi(-1)$. Substituting this into \eqref{eqn-prop-2} we can write
\begin{align}
B'&=-\varphi(-1)\psi_3(4)\overline{\psi}_6(2)\sum_{a=0}^{p-2}\pi^{(p-1)S}
~\Gamma_p\left(\left\langle\frac{5}{6}-\frac{a}{p-1}\right\rangle\right)^2\Gamma_p\left(\left\langle\frac{2}{3}+\frac{a}{p-1}\right\rangle\right)^4\notag\\
&\times\Gamma_p\left(\left\langle\frac{1}{12}-\frac{a}{p-1}\right\rangle\right)
\Gamma_p\left(\left\langle\frac{7}{12}-\frac{a}{p-1}\right\rangle\right).
\end{align}
We have \begin{align}\label{S-1}
S&=2\left\langle\frac{5}{6}-\frac{a}{p-1}\right\rangle+4\left\langle\frac{2}{3}+\frac{a}{p-1}\right\rangle
+\left\langle\frac{1}{6}-\frac{2a}{p-1}\right\rangle-\frac{1}{2}\notag\\
&=4-2\left\lfloor\frac{5}{6}-\frac{a}{p-1}\right\rfloor-\left\lfloor\frac{1}{6}-\frac{2a}{p-1}\right\rfloor
-4\left\lfloor\frac{2}{3}+\frac{a}{p-1}\right\rfloor.
\end{align}
Also, we can easily see that
$\left\lfloor\frac{1}{6}-\frac{2a}{p-1}\right\rfloor=\left\lfloor\frac{1}{12}-\frac{a}{p-1}\right\rfloor
+\left\lfloor\frac{7}{12}-\frac{a}{p-1}\right\rfloor$. Then we have 
$S=4-2\left\lfloor\frac{5}{6}-\frac{a}{p-1}\right\rfloor-\left\lfloor\frac{1}{12}-\frac{a}{p-1}\right\rfloor
-\left\lfloor\frac{7}{12}-\frac{a}{p-1}\right\rfloor
-4\left\lfloor\left\langle-\frac{1}{3}\right\rangle+\frac{a}{p-1}\right\rfloor.$ This gives
\begin{align}
B'&=-p^4\varphi(-1)\psi_3(4)\overline{\psi}_6(2)\sum_{a=0}^{p-2}
(-p)^{-2\left\lfloor\frac{5}{6}-\frac{a}{p-1}\right\rfloor-\left\lfloor\frac{1}{12}-\frac{a}{p-1}\right\rfloor
-\left\lfloor\frac{7}{12}-\frac{a}{p-1}\right\rfloor
-4\left\lfloor\left\langle-\frac{1}{3}\right\rangle+\frac{a}{p-1}\right\rfloor}\notag\\
&\times\Gamma_p\left(\left\langle\frac{5}{6}-\frac{a}{p-1}\right\rangle\right)^2\Gamma_p\left(\left\langle\frac{2}{3}+\frac{a}{p-1}\right\rangle\right)^4\notag\\
&\times\Gamma_p\left(\left\langle\frac{1}{12}-\frac{a}{p-1}\right\rangle\right)
\Gamma_p\left(\left\langle\frac{7}{12}-\frac{a}{p-1}\right\rangle\right)
\notag\\
&=p^4(p-1)\varphi(-1)\psi_3(4)\overline{\psi}_6(2)~
{_4G_4}\left[\begin{array}{cccc}
             \frac{5}{6}, & \frac{5}{6}, & \frac{1}{12}, & \frac{7}{12}\vspace{.1 cm } \\
             \frac{1}{3}, & \frac{1}{3}, & \frac{1}{3}, & \frac{1}{3}
           \end{array}|1\right].\notag
\end{align}
This completes the proof of the proposition.
\end{proof}
\begin{proposition}\label{prop-mt-2}
Let $p>3$ be an odd prime such that $p\not\equiv1\pmod{3}$. Then 
\begin{align*}
&\frac{\varphi(-1)}{p(1-p)}\sum_{a=0}^{p-2}\frac{g(\varphi\omega^{a})^2g(\overline{\omega}^a)^4g(\varphi\omega^{2a})}{g(\varphi)}~\overline{\omega}^a(4)\\
&={_{12}G_{12}}\left[\begin{array}{cccccccccccc}
             0, & \frac{1}{3}, & \frac{2}{3}, & 0, & \frac{1}{3}, & \frac{2}{3}, & 0, & \frac{1}{3}, & \frac{2}{3}, &
             0, & \frac{1}{3}, & \frac{2}{3} \vspace{.1 cm } \\
             \frac{1}{6}, & \frac{1}{2}, & \frac{5}{6}, & \frac{1}{12}, & \frac{1}{6}, & \frac{1}{4},
             & \frac{5}{12}, & \frac{1}{2}, & \frac{7}{12}, & \frac{3}{4}, & \frac{5}{6}, & \frac{11}{12}
           \end{array}\mid1\right]
\end{align*}
\end{proposition}
\begin{proof}
Let $D=\displaystyle\sum_{a=0}^{p-2}\frac{g(\varphi\omega^{a})^2g(\overline{\omega}^a)^4g(\varphi\omega^{2a})}{g(\varphi)}~\overline{\omega}^a(4)$. Since $(p-1,3)=1$ so taking the automorphism $a\rightarrow 3a$ we have
\begin{align}
D=\sum_{a=0}^{p-2}\frac{g(\varphi\omega^{3a})^2g(\overline{\omega}^{3a})^4g(\varphi\omega^{6a})}{g(\varphi)}~\overline{\omega}^{3a}(4).\notag
\end{align}
Applying Davenport Hasse relation for $m=2$ and $\psi=\omega^{3a}, \omega^{6a}$ successively, we have
\begin{align}
g(\varphi\omega^{3a})&=\frac{g(\omega^{6a})\omega^{3a}(2^{-2})g(\varphi)}{g(\omega^{3a})},\notag\\
g(\varphi\omega^{6a})&=\frac{g(\omega^{12a})\omega^{6a}(2^{-2})g(\varphi)}{g(\omega^{6a})}.\notag
\end{align}
Substituting these two expressions into the above expression of $D$ we can write
\begin{align}\label{eqn-prop-3}
D&=\sum_{a=0}^{p-2}\frac{g(\omega^{6a})^2\omega^{6a}(2^{-2})g(\varphi)^2g(\overline{\omega}^{3a})^4}{g(\omega^{3a})^2}\frac{g(\omega^{12a})\omega^{6a}(2^{-2})}{g(\omega^{6a})}~\overline{\omega}^{3a}(4).
\end{align}
Lemma \ref{lemma-1} gives $g(\varphi)^2=p\varphi(-1)$. Hence, using this in \eqref{eqn-prop-3} we obtain
\begin{align}
D&=p\varphi(-1)\sum_{a=0}^{p-2}\frac{g(\overline{\omega}^{3a})^4g(\omega^{6a})
g(\omega^{12a})}{g(\omega^{3a})^2}~\overline{\omega}^{30a}(2).
\end{align}
Now, using Gross Koblitz formula we have 
\begin{align}
D&=p\varphi(-1)\sum_{a=0}^{p-2}\overline{\omega}^{30a}(2)\pi^{(p-1)T}~
\frac{\Gamma_p(\langle\frac{3a}{p-1}\rangle)^4\Gamma_p(\langle-\frac{6a}{p-1}\rangle)
\Gamma_p(\langle-\frac{12a}{p-1}\rangle)}{\Gamma_p(\langle-\frac{3a}{p-1}\rangle)^2},
\end{align}
where $T=4\langle\frac{3a}{p-1}\rangle+\langle\frac{-6a}{p-1}\rangle+\langle\frac{-12a}{p-1}\rangle
-2\langle\frac{-3a}{p-1}\rangle$.
Then applying Lemma \ref{lemma-5} and simplifying we obtain
\begin{align}\label{eqn-prop-4}
D&=p\varphi(-1)\sum_{a=0}^{p-2}\pi^{(p-1)T}~~\frac{\Gamma_p(\langle\frac{a}{p-1}\rangle)^4
\Gamma_p(\langle\frac{1}{3}+\frac{a}{p-1}\rangle)^4\Gamma_p(\langle\frac{2}{3}+\frac{a}{p-1}\rangle)^4}
{\Gamma_p(\frac{1}{3})^4\Gamma_p(\frac{2}{3})^4}\notag\\
&\times\frac{\Gamma_p(\langle\frac{1}{6}-\frac{a}{p-1}\rangle)\Gamma_p(\langle\frac{1}{2}-\frac{a}{p-1}\rangle)
\Gamma_p(\langle\frac{5}{6}-\frac{a}{p-1}\rangle)\Gamma_p(\langle\frac{1}{12}-\frac{a}{p-1}\rangle)}
{\Gamma_p(\frac{1}{6})\Gamma_p(\frac{1}{2})\Gamma_p(\frac{5}{6})\Gamma_p(\frac{1}{12})}\notag\\
&\times\frac{\Gamma_p(\langle\frac{1}{6}-\frac{a}{p-1}\rangle)\Gamma_p(\langle\frac{1}{4}-\frac{a}{p-1}\rangle)
\Gamma_p(\langle\frac{5}{12}-\frac{a}{p-1}\rangle)\Gamma_p(\langle\frac{1}{2}-\frac{a}{p-1}\rangle)
\Gamma_p(\langle\frac{7}{12}-\frac{a}{p-1}\rangle)}
{\Gamma_p(\frac{1}{6})\Gamma_p(\frac{1}{4})\Gamma_p(\frac{5}{12})\Gamma_p(\frac{1}{2})\Gamma_p(\frac{7}{12})}\notag\\
&\times\frac{\Gamma_p(\langle\frac{3}{4}-\frac{a}{p-1}\rangle)
\Gamma_p(\langle\frac{5}{6}-\frac{a}{p-1}\rangle)\Gamma_p(\langle\frac{11}{12}-\frac{a}{p-1}\rangle)}
{\Gamma_p(\frac{3}{4})\Gamma_p(\frac{5}{6})\Gamma_p(\frac{11}{12})}.
\end{align}
We have \begin{align}
T&=4\left\langle\frac{3a}{p-1}\right\rangle+\left\langle\frac{-6a}{p-1}\right\rangle
+\left\langle\frac{-12a}{p-1}\right\rangle
-2\left\langle\frac{-3a}{p-1}\right\rangle\notag\\
&=-4\left\lfloor\frac{3a}{p-1}\right\rfloor-\left\lfloor\frac{-6a}{p-1}\right\rfloor
-\left\lfloor\frac{-12a}{p-1}\right\rfloor+2\left\lfloor\frac{-3a}{p-1}\right\rfloor.
\end{align}
Using Lemma \ref{D1} and Lemma \ref{D2} we obtain
\begin{align}\label{eqn-prop-5}
T&=-4\left\lfloor\frac{a}{p-1}\right\rfloor-4\left\lfloor\frac{1}{3}+\frac{a}{p-1}\right\rfloor
-4\left\lfloor\frac{2}{3}+\frac{a}{p-1}\right\rfloor-2\left\lfloor\frac{1}{6}-\frac{a}{p-1}\right\rfloor\notag\\
&-2\left\lfloor\frac{1}{2}-\frac{a}{p-1}\right\rfloor-2\left\lfloor\frac{5}{6}-\frac{a}{p-1}\right\rfloor
-\left\lfloor\frac{1}{12}-\frac{a}{p-1}\right\rfloor-\left\lfloor\frac{1}{4}-\frac{a}{p-1}\right\rfloor-
\left\lfloor\frac{5}{12}-\frac{a}{p-1}\right\rfloor\notag\\
&-\left\lfloor\frac{7}{12}-\frac{a}{p-1}\right\rfloor-\left\lfloor\frac{3}{4}-\frac{a}{p-1}\right\rfloor-
\left\lfloor\frac{11}{12}+\frac{a}{p-1}\right\rfloor.
\end{align}
Finally, using \eqref{eqn-prop-5} into \eqref{eqn-prop-4} we obtain
\begin{align}
&\frac{\varphi(-1)D}{p(p-1)}\notag\\&=
-{_{12}G_{12}}\left[\begin{array}{cccccccccccc}
             0, & \frac{1}{3}, & \frac{2}{3}, & 0, & \frac{1}{3}, & \frac{2}{3}, & 0, & \frac{1}{3}, & \frac{2}{3}, &
             0, & \frac{1}{3}, & \frac{2}{3} \vspace{.1 cm } \\
             \frac{1}{6}, & \frac{1}{2}, & \frac{5}{6}, & \frac{1}{12}, & \frac{1}{6}, & \frac{1}{4},
             & \frac{5}{12}, & \frac{1}{2}, & \frac{7}{12}, & \frac{3}{4}, & \frac{5}{6}, & \frac{11}{12}
           \end{array}\mid1\right].\notag
\end{align}
This completes the proof of the proposition.
\end{proof}
\subsection{Reductions of Twisted Kloosterman sums}
The following proposition expresses the sum $S(m,\varphi)$ as a symmetric sum of quadratic characters.
\begin{proposition}\cite[Proposition 3.1]{DGP}\label{proposition-1} For an integer $m\geq1,$ we have
\begin{align*}
S(m+1,\varphi)=p\varphi(-1)\sum_{x_1,\ldots,x_m\in\mathbb{F}_p^{\times}}\varphi(x_1+x_2+\cdots+x_m+1)\varphi(\overline{x_1}+\overline{x_2}+
\cdots+\overline{x_m}+1)
\end{align*}
\end{proposition}
Let $F:\mathbb{F}_p^{\times}\rightarrow\mathbb{C}$ be a function defined by
$$F(a):=\sum_{x,y\in\mathbb{F}_p^{\times}}\varphi(x+y+a+1)\varphi(\overline{x}+\overline{y}
+\overline{a}+1).$$
The following two lemmas build connections between the function $F$, and special values of the hypergeometric functions ${_2F_1}(x).$ 
\begin{lemma}\label{lemma-8}\cite[Lemma 3.2]{DGP}
For $a\neq0,\pm 1$, we have
\begin{align*}
F(a)=p^2\varphi(a)~{_2F_1}(-a)^2.
\end{align*}
\end{lemma}
\begin{lemma}\label{lemma-9}\cite[Lemma 3.3]{DGP}
We have that 
\begin{align*}
\varphi(-1)F(1)&=p^2\varphi(-1)~{_2F_1}(-1)^2-p,\\
\varphi(-1)F(-1)&=p^2~{_2F_1}(1)^2-p.
\end{align*}
\end{lemma}
\section{Proof of main results}
In this section, we are ready to prove our main results. We first proof Theorem \ref{MT-1}. To prove this theorem we express the 4-th quadratic Kloosterman sheaf sum in terms of Gauss sums and then we apply Proposition \ref{prop-mt-1}.
\begin{proof}[Proof of Theorem \ref{MT-1}]
By definition we have 
\begin{align*}
T_{4,\varphi}=\sum_{a\in\mathbb{F}_p^{\times}}\varphi(a)(g(a)^4+g(a)^3h(a)+g(a)^2h(a)^2+g(a)h(a)^3+h(a)^4).
\end{align*}
Then using the relations $g(a)+h(a)=-K(a)$ and $g(a)h(a)=p$ we can write
$g(a)^4+g(a)^3h(a)+g(a)^2h(a)^2+g(a)h(a)^3+h(a)^4=K(a)^4-3pK(a)^2+p^2$. This gives
\begin{align}\label{eqn-3}
T_{4,\varphi}&=\sum_{a\in\mathbb{F}_p^{\times}}\varphi(a)(K(a)^4-3pK(a)^2+p^2)\notag\\
&=\sum_{a\in\mathbb{F}_p^{\times}}\varphi(a)K(a)^4-3p\sum_{a\in\mathbb{F}_p^{\times}}\varphi(a)K(a)^2+
p^2\sum_{a\in\mathbb{F}_p^{\times}}\varphi(a)\notag\\
&=S(4,\varphi)-3pS(2,\varphi).
\end{align}
Here note that the sum $p^2\displaystyle\sum_{a\in\mathbb{F}_p^{\times}}\varphi(a)=0$ by first part of Lemma \ref{lemma-3}. Applying Proposition \ref{proposition-1} with $m=3,1$ successively, we obtain
\begin{align}\label{eqn-4}
S(4,\varphi)&=p\varphi(-1)\sum_{x,y,z\in\mathbb{F}_p^{\times}}\varphi(x+y+z+1)
\varphi(\overline{x}+\overline{y}+\overline{z}+1)\notag\\
&=p\varphi(-1)\sum_{x\in\mathbb{F}_p^{\times}}F(x)
\end{align}
and 
\begin{align}\label{eqn-5}
S(2,\varphi)&=p\varphi(-1)\sum_{x\in\mathbb{F}_p^{\times}}\varphi(x+1)\varphi(\overline{x}+1)\notag\\
&=p\varphi(-1)\sum_{x\in\mathbb{F}_p^{\times}, x\neq-1}\varphi(x)=-p.
\end{align}
Using \eqref{eqn-4} and \eqref{eqn-5} in \eqref{eqn-3} we obtain
\begin{align*}
T_{4,\varphi}=p\varphi(-1)\sum_{x\in\mathbb{F}_p^{\times}}F(x)+3p^2.
\end{align*}
Now applying Lemma \ref{lemma-8} and Lemma \ref{lemma-9} in the above expression of $T_{4,\varphi}$ we have
\begin{align*}
\frac{T_{4,\varphi}}{p}&=p+p^2\sum_{x\in\mathbb{F}_p^{\times}}\varphi(x)~{_2F_1}(x)^2\\
&=p+p^2\sum_{x\in\mathbb{F}_p}\varphi(x)~{_2F_1}(x)^2.
\end{align*}
Replacing $x$ by $\frac{1-t}{2}$ we have
\begin{align}\label{eqn-6}
\frac{T_{4,\varphi}}{p}&=p+p^2\sum_{t\in\mathbb{F}_p}\varphi\left(\frac{1-t}{2}\right)~{_2F_1}\left(\frac{1-t}{2}\right)^2\notag\\
&=p+p^2\varphi(2)~{_2F_1}\left(\frac{1}{2}\right)^2+p^2~{_2F_1}(1)^2+
p^2\sum_{t\in\mathbb{F}_p^{\times}, t\neq\pm1}\varphi\left(\frac{1-t}{2}\right)~{_2F_1}\left(\frac{1-t}{2}\right)^2.
\end{align}
Let $A=\displaystyle\varphi(2)\sum_{t\in\mathbb{F}_p^{\times}, t\neq\pm1}\varphi(1-t)~{_2F_1}\left(\frac{1-t}{2}\right)^2.$ Then using Lemma \ref{lemma 10} we have
\begin{align}
A&=\varphi(2)\sum_{t\in\mathbb{F}_p^{\times}, t\neq\pm1}\varphi(1-t)
\left(\frac{1}{p}+\varphi(t^2-1){_3F_2}\left(\frac{1}{1-t^2}\right)\right)\notag\\
&=\frac{\varphi(2)}{p}\sum_{t\in\mathbb{F}_p^{\times}, t\neq\pm1}\varphi(1-t)+\varphi(2)
\sum_{t\in\mathbb{F}_p^{\times}, t\neq\pm1}\varphi(1-t)\varphi(t^2-1){_3F_2}\left(\frac{1}{1-t^2}\right)\notag\\
&=-\frac{1}{p}-\frac{\varphi(2)}{p}+\frac{\varphi(2)}{p}\sum_{t\in\mathbb{F}_p}\varphi(1-t)
+\varphi(-2)\sum_{t\in\mathbb{F}_p^{\times}, t\neq\pm1}\varphi(1+t){_3F_2}\left(\frac{1}{1-t^2}\right).\notag
\end{align}
In the above expression the second sum vanishes by first part of Lemma \ref{lemma-3}. Therefore, we have 
\begin{align}
A&=-\frac{1}{p}-\frac{\varphi(2)}{p}+\varphi(-2)\sum_{t\in\mathbb{F}_p^{\times}, t\neq\pm1}\varphi(1+t)
\sum_{\chi\in\widehat{\mathbb{F}_p^{\times}}}{\varphi\chi\choose\chi}^3\chi\left(\frac{1}{1-t^2}\right)\notag\\
&=-\frac{1}{p}-\frac{\varphi(2)}{p}+\varphi(-2)\sum_{t\in\mathbb{F}_p^{\times}, t\neq\pm1}\varphi(1+t)
\sum_{\chi\in\widehat{\mathbb{F}_p^{\times}}}{\varphi\chi\choose\chi}^3\overline{\chi}(1-t^2)\notag.
\end{align}
Then adjusting the term under summation for $t=0,\pm1$ we can write
\begin{align}\label{eqn-7}
A&=-\frac{1}{p}-\frac{\varphi(2)}{p}+\varphi(-2)\sum_{t\in\mathbb{F}_p}\varphi(1+t)
\sum_{\chi\in\widehat{\mathbb{F}_p^{\times}}}{\varphi\chi\choose\chi}^3\overline{\chi}(1-t^2)
-\varphi(-2)\sum_{\chi\in\widehat{\mathbb{F}_p^{\times}}}{\varphi\chi\choose\chi}^3\notag\\
&=-\frac{1}{p}-\frac{\varphi(2)}{p}-\varphi(-2){_3F_2}(1)+\varphi(-2)\sum_{t\in\mathbb{F}_p}\varphi(1+t)
\sum_{\chi\in\widehat{\mathbb{F}_p^{\times}}}{\varphi\chi\choose\chi}^3\overline{\chi}(1-t^2).
\end{align}
Now, consider the sum
\begin{align}\label{eqn-8}
B&=\sum_{t\in\mathbb{F}_p}\varphi(1+t)
\sum_{\chi\in\widehat{\mathbb{F}_p^{\times}}}{\varphi\chi\choose\chi}^3\overline{\chi}(1-t^2)\notag\\
&=\sum_{\chi\in\widehat{\mathbb{F}_p^{\times}}}{\varphi\chi\choose\chi}^3\sum_{t\in\mathbb{F}_p}
\varphi \overline{\chi}(1+t)\overline{\chi}(1-t).
\end{align}
If we use \eqref{eqn-J} then we obtain $\displaystyle\sum_{t\in\mathbb{F}_p}
\varphi \overline{\chi}(1+t)\overline{\chi}(1-t)=\varphi\overline{\chi}^2(2)J(\varphi\overline{\chi},\overline{\chi})$. Now, using this in \eqref{eqn-8}  and then applying \eqref{binomial} we obtain
\begin{align}
B&=\varphi(2)\sum_{\chi\in\widehat{\mathbb{F}_p^{\times}}}{\varphi\chi\choose\chi}^3\overline{\chi}^2(2)
J(\varphi\overline{\chi},\overline{\chi})\notag\\
&=p\varphi(2)\sum_{\chi\in\widehat{\mathbb{F}_p^{\times}}}{\varphi\chi\choose\chi}^3
{\varphi\overline{\chi}\choose\chi}\chi\left(-\frac{1}{4}\right).\notag
\end{align}
Now, \eqref{property-1} implies ${\varphi\overline{\chi}\choose\chi}=\chi(-1){\varphi\chi^2\choose\chi}$. Using this relation in the above sum we have
\begin{align}
B&=p\varphi(2)\sum_{\chi\in\widehat{\mathbb{F}_p^{\times}}}{\varphi\chi\choose\chi}^3{\varphi\chi^2\choose\chi}
\chi\left(\frac{1}{4}\right).\notag
\end{align}
Using \eqref{property-2} and the fact that $g(\varepsilon)=-1$ we deduce that 
\begin{align}
B&=\frac{\varphi(2)}{p^3}\sum_{\chi\in\widehat{\mathbb{F}_p^{\times}}}
\frac{g(\varphi\chi)^2g(\overline{\chi})^4g(\varphi\chi^2)}{g(\varphi)^3}\overline{\chi}(4)-\frac{p-1}{p^3}\varphi(-2).\notag
\end{align}
Using Lemma \ref{lemma-1} we have $g(\varphi)^2=p\varphi(-1)$. Therefore, using this and then 
taking $\chi=\omega^a$ we can write
\begin{align}\label{eqn-9}
B=\frac{\varphi(-2)}{p^4}\sum_{a=0}^{p-2}
\frac{g(\varphi\omega^a)^2g(\overline{\omega}^a)^4g(\varphi\omega^{2a})}{g(\varphi)}~\overline{\omega}^a(4)-\frac{p-1}{p^3}\varphi(-2).
\end{align}
Now, substituting \eqref{eqn-9} into \eqref{eqn-7} and then again substituting the resultant identity of $A$ into 
\eqref{eqn-6} we deduce that 
\begin{align}\label{eqn-10}
\frac{T_{4,\varphi}}{p}&=-\frac{p-1}{p}-p\varphi(2)+p^2\varphi(2){_2F_1}\left(\frac{1}{2}\right)^2+p^2{_2F_1}(1)^2
-p^2\varphi(-2){_3F_2}(1)\notag\\
&+\frac{1}{p^2}\sum_{a=0}^{p-2}
\frac{g(\varphi\omega^a)^2g(\overline{\omega}^a)^4g(\varphi\omega^{2a})}{g(\varphi)}~\overline{\omega}^a(4).
\end{align}
Now, if $p\equiv1\pmod{12}$, then by Proposition \ref{value-1} we have 
\begin{align}\label{eqn-11}
{_2F_1}\left(\frac{1}{2}\right)=\frac{2x(-1)^{\frac{x+y+1}{2}}}{p},
\end{align}
and by Proposition \ref{value-2} we have
\begin{align}\label{eqn-13}
{_3F_2}(1)=\frac{4x^2-2p}{p^2},
\end{align} where $p=x^2+y^2$ with $x\equiv1\pmod{2}$,
and Proposition \ref{value-3} gives
\begin{align}\label{eqn-12}
{_2F_1(1)}=-\frac{\varphi(-1)}{p}.
\end{align}
Now, substituting \eqref{eqn-11}, \eqref{eqn-13} and \eqref{eqn-12} into \eqref{eqn-10} we obtain
\begin{align}\label{eqn-14}
\frac{T_{4,\varphi}}{p}&=\frac{1}{p}-p\varphi(2)+4x^2\varphi(2)-\varphi(-2)(4x^2-2p)\notag\\
&+\frac{1}{p^2}\sum_{a=0}^{p-2}
\frac{g(\varphi\omega^a)^2g(\overline{\omega}^a)^4g(\varphi\omega^{2a})}{g(\varphi)}~\overline{\omega}^a(4).
\end{align}
Applying Proposition \ref{prop-mt-1} in the last summation involve in \eqref{eqn-14}  we prove \eqref{case-1}.
Again, if $p\equiv7\pmod{12}$, then by Proposition \ref{value-1} we have 
\begin{align}\label{eqn-15}
{_2F_1}\left(\frac{1}{2}\right)=0,
\end{align}
and by Proposition \ref{value-2} we have
\begin{align}\label{eqn-16}
{_3F_2}(1)=0.
\end{align}
Then substituting \eqref{eqn-12}, \eqref{eqn-15} and \eqref{eqn-16} into \eqref{eqn-10} we obtain
\begin{align}\label{eqn-17}
\frac{T_{4,\varphi}}{p}&=\frac{1}{p}-p\varphi(2)+\frac{1}{p^2}\sum_{a=0}^{p-2}
\frac{g(\varphi\omega^a)^2g(\overline{\omega}^a)^4g(\varphi\omega^{2a})}{g(\varphi)}~\overline{\omega}^a(4).
\end{align}
Again, applying Proposition \ref{prop-mt-1} in the last summation involve in \eqref{eqn-17}  we prove \eqref{case-2}. This completes the proof of the theorem.
\end{proof}
\begin{proof}[Proof of Theorem \ref{MT-2}]
From \eqref{eqn-10}  we write
\begin{align}\label{eqn-18}
\frac{T_{4,\varphi}}{p}&=-\frac{p-1}{p}-p\varphi(2)+p^2\varphi(2){_2F_1}\left(\frac{1}{2}\right)^2+p^2{_2F_1}(1)^2
-p^2\varphi(-2){_3F_2}(1)\notag\\
&+\frac{1}{p^2}\sum_{a=0}^{p-2}
\frac{g(\varphi\omega^a)^2g(\overline{\omega}^a)^4g(\varphi\omega^{2a})}{g(\varphi)}~\overline{\omega}^a(4).\end{align}
Since $p\equiv5\pmod{12}$, then by Proposition \ref{value-1} we have 
\begin{align}\label{eqn-19}
{_2F_1}\left(\frac{1}{2}\right)=\frac{2x(-1)^{\frac{x+y+1}{2}}}{p},
\end{align}
and by Proposition \ref{value-2} we have
\begin{align}\label{eqn-20}
{_3F_2}(1)=\frac{4x^2-2p}{p^2},
\end{align}
where $p=x^2+y^2$ with $x\equiv1\pmod{2}$.
Then substituting \eqref{eqn-19}, \eqref{eqn-20}, and \eqref{eqn-12} into \eqref{eqn-18} we have
\begin{align}\label{eqn-21}
\frac{T_{4,\varphi}}{p}&=\frac{1}{p}-p\varphi(2)+4x^2\varphi(2)-\varphi(-2)(4x^2-2p)\notag\\
&+\frac{1}{p^2}\sum_{a=0}^{p-2}
\frac{g(\varphi\omega^a)^2g(\overline{\omega}^a)^4g(\varphi\omega^{2a})}{g(\varphi)}~\overline{\omega}^a(4).
\end{align}
Now, applying Proposition \ref{prop-mt-2} in the last summation involve in \eqref{eqn-21}  we prove \eqref{case-3}.
Also, if $p\equiv11\pmod{12}$, then we have ${_2F_1}\left(\frac{1}{2}\right)={_3F_2}(1)=0$ and ${_2F_1(1)}=-\frac{\varphi(-1)}{p}$.  Hence putting these values into \eqref{eqn-18} and again applying
 Proposition \ref{prop-mt-2} we obtain \eqref{case-4}. This completes the proof of the theorem.
\end{proof}
\begin{proof}[Proof of Theorem \ref{MT-3} and \ref{MT-4}] 
In \cite{ahl-ono-1} and \cite{ahl-ono-2}, Ahlgren and Ono proved 
\begin{align}\label{eq-1}
p^3~{_4F_3}(1)=-a(p)-p.
\end{align}
From \cite[Theorem 1.1]{DGP} we have
\begin{align}\label{eq-2}
\frac{T_4}{p}=p^3~{_4F_3}(1)+p.
\end{align}
Combining \eqref{eq-1} and \eqref{eq-2} we obtain
\begin{align}\label{eq-3}
\frac{T_4}{p}=-a(p).
\end{align}
Since $f_2$ is the quadratic twist of $f_1$, so $b(p)=\varphi(-1)a(p)$. This gives
\begin{align}\label{eq-4}
\frac{T_4}{p}=-\varphi(-1)b(p).
\end{align}
Then using Theorem \ref{MT-1} in \eqref{eq-3} and \eqref{eq-4} we prove Theorem \ref{MT-3}. Finally, applying Theorem \ref{MT-2} in \eqref{eq-3} and \eqref{eq-4} we prove Theorem \ref{MT-4}.
\end{proof}

\end{document}